\documentclass{article}

\usepackage{amsmath}

\renewcommand{\eqref}[1]{(\ref{#1})}

\usepackage{amsfonts}
\usepackage{graphicx}

\newcommand{\Laplacian}{\mathrm{\Delta}}

\newcommand{\pder}[2]{\frac{\partial{#1}}{\partial{#2}}}

\renewcommand{\div}{\nabla_x\cdot}
\newcommand{\grad}{\nabla_x}

\newcommand{\figuresdir}{figurearxiv}

\title{High order relaxed schemes for nonlinear reaction diffusion
  problems in nonconservative form}

\author{F. Cavalli and M. Semplice$^*$\\
Dipartimento di Matematica, Universit\`a di Milano\\
Via Saldini 50, I-20123 Milano, Italy\\
$^*$E-mail: semplice@mat.unimi.it}

\begin{document}
\maketitle

\begin{abstract}
Different relaxation approximations to partial differential equations,
including conservation laws, Hamilton-Jacobi equations,
convection-diffusion problems, gas dynamics problems, have been
recently proposed. The present paper focuses onto diffusive relaxed
schemes for the numerical approximation of nonlinear reaction
diffusion equations. We choose here a nonstandard relaxation scheme
that allow the treatment of diffusion equations in their
nonconservative form. A comparison with the traditional approach in
the case of conservative equations is also included.
  High order methods are obtained by coupling ENO
and WENO schemes for space discretization with IMEX schemes for time
integration, where the implicit part can be explicitly solved at a
linear cost.  To illustrate the high accuracy and good properties of
the proposed numerical schemes, also in the degenerate case, we
consider various examples in one and two dimensions: the
Fisher-Kolmogoroff equation, the porous-Fisher equation and the porous
medium equation with strong absorption. Moreover we show a test on a
system of PDEs that describe an ecological model for the dispersal and
settling of animal populations.
\end{abstract}


\section{Introduction}

The main purpose of this work is to approximate solutions of a nonlinear,
possibly degenerate, reaction-diffusion equation of the form
\begin{equation}  \label{eq:degparab}
\frac{\partial u}{\partial t} = D \nabla\cdot(A(u) \nabla u) + g(u)
\end{equation}
for $ x\in \Omega \subseteq\mathbb{R}^d,~~d\geq 1,~~ t\geq0$, with
initial condition, $u(x,0)=u_0(x)$ and with suitable boundary
conditions,  
where the function $A$ is a non-negative bounded 
function defined on $\mathbb{R}$. The equation is degenerate if
$A(0)=0$ and in this case the solutions often become non-smooth in
finite time, developing fronts and discontinuities \cite{Aro70}. 
Finally, the coefficient $D$ is a diffusivity
coefficient and the function $g(u)$ is the reaction term.

Equations like \eqref{eq:degparab} are relevant in describing
biological and physical processes and are involved in the modelling of
population growth and dispersal \cite{GSM00,TLB06}, of waves of
concentration of chemical substances in living organisms, of the
motion of viscous fluids (see e.g. \cite{Murray1}). Also, systems of
equations like \eqref{eq:degparab}, coupled via the reaction terms,
are capable of modelling the cyclic Belousov-Zhabotinskii reactions
and the pattern formation on the wings of butterflies and the coat of
mammals (see \cite{Murray2} and references therein).

This paper is organized as follows. In section
\ref{section:relaxation} we introduce the relaxation approximation of
nonlinear diffusion problems, in section \ref{sec:scheme} we describe
the fully discrete relaxed numerical scheme in the reaction-diffusion
case. In section \ref{sec:tests} we report several numerical tests,
in one and two space dimensions. The last application shown is about a
relevant ecological model.

\section{Relaxation approximation of nonlinear diffusion}
\label{section:relaxation}

The schemes proposed in the present work are a generalization of those
proposed and studied in \cite{NP98}, and fit in the wider class of the
schemes associated to relaxation approximations of partial differential
equations, started after the well-known case studied in \cite{JX95}
for hyperbolic conservation laws.
In the case of the nonlinear
diffusion operator in \eqref{eq:degparab}, an additional variable
$\vec{v}(x,t)\in\mathbb{R}^d$ and a positive parameter $\varepsilon$
are introduced and the following relaxation system is considered:
\begin{equation}\label{eq:sysrel}
\left\{
\begin{array}{ll}
\displaystyle
\frac{\partial u}{\partial t} + \mathrm{div}(\vec{v}) = g(u)\\
\\
\displaystyle
\frac{\partial \vec{v}}{\partial t} +\frac{D}{\varepsilon} A(u) \nabla u   = 
    -\frac 1\varepsilon \vec{v} 
\end{array}
\right.
\end{equation}
Now, formally, in the small relaxation limit, $\varepsilon \rightarrow 
0^+$, system \eqref{eq:sysrel} approximates to leading order
equation \eqref{eq:degparab}.

We point out that in this relaxation system we do not need a third
equation to ``relax'' the non-linearity of $A$, as was the case in
\cite{CNPSdegdiff}. 

In the previous system the parameter $\varepsilon$ has physical
dimensions of time and represents the so-called relaxation time.
Furthermore, each component of $\vec{v}$ has the dimension of $u$
times a velocity and $\varphi$ is a velocity.  The inverse of
$\varepsilon$ gives the rate at which $v$ decays onto $-A(u)\nabla u$
in the evolution of the variable $\vec{v}$ governed by the stiff
second equation of \eqref{eq:sysrel}.

In order to justify the convergence of the proposed relaxation scheme
to equation \eqref{eq:degparab}, we study the Chapman-Enskog expansion
of the solutions of \eqref{eq:sysrel}. We consider for simplicity the
case with only one spatial variable and temporarily assume that $A$ is
sufficiently regular. From the second equation of
\eqref{eq:sysrel} we obtain
\[ v=-A(u)u_x -\varepsilon v_t
    = -A(u)u_x -\varepsilon (A(u)u_x)_t + O(\varepsilon^2)
\]
Using the relation above and \eqref{eq:sysrel} again, one computes
\[\begin{split}
(A(u)u_x)_t = A'(u)u_tu_x + A(u)u_{tx}
=A'(u)(-v_x)u_x + A(u)(-v_x)_x\\
=-(A(u)v_x)_x=\left(A(u)\left(A(u)u_x\right)_x\right)_x+O(\varepsilon)
\end{split}\]
Finally we substitute in the first equation of \eqref{eq:sysrel},
obtaining
\[
u_t-\left(A(u)u_x\right)_x = \varepsilon
\left(A(u)\left(A(u)u_x\right)_x\right)_{xx} + O(\varepsilon^2)
\]

The idea of employing a semilinear hyperbolic relaxation system to
derive numerical schemes for (degenerate) diffusion equation was first
introduced in \cite{NP00} for the purely diffusive case in the
conservative form $u_t=\Laplacian(p(u))$. High order schemes based on
such approach have been studied in \cite{CNPSdegdiff}. As in that paper
we introduce a parameter $\varphi$ which allows to move the
stiff terms $\frac{D}{\varepsilon}\nabla p(u)$ to the right hand side,
without losing the hyperbolicity of the system:
\begin{equation}\label{eq:sysrelphi}
\left\{
\begin{array}{ll}
\displaystyle
\frac{\partial u}{\partial t} + \mathrm{div}(\vec{v}) = g(u)\\
\\
\displaystyle
\frac{\partial \vec{v}}{\partial t} + \varphi^2 \nabla u = 
    -\frac 1\varepsilon \vec{v} + 
    \left(\varphi^2 \nabla u - \frac{D}{\varepsilon} A(u) \nabla u \right)
\end{array}
\right.
\end{equation}
Now, the characteristic velocities of the hyperbolic
left hand side are given by $\pm\varphi$ and are no longer stiff.

We point out that degenerate parabolic equations often model physical
situations where free boundaries and discontinuities are relevant: we
expect that schemes for hyperbolic systems will be able to reproduce
faithfully these details of the solution.  One of the main properties
of \eqref{eq:sysrel} consists in the semilinearity of the system, that is
all the nonlinear terms are in the (stiff) source terms, while the
differential operator is linear. Hence, the solution of the convective
part requires neither Riemann solvers nor the computation of the
characteristic structure at each time step, since the eigenstructure
of the system is constant in time.  Moreover, the relaxation
approximation does not exploit the form of the nonlinear function $A$
and hence it gives rise to a numerical scheme that, to a large extent,
is independent of it, resulting in a very versatile tool.

We also anticipate here that, in the relaxed case
(i.e. $\varepsilon=0$), the stiff source terms can be integrated
solving a system that is already in triangular form and then it does
not require iterative solvers.

Obviously, when \eqref{eq:degparab} can be put in the form
$u_t+\Laplacian(p(u))=g(u)$, the relaxation approximation of the
parabolic part described in \cite{CNPSdegdiff} can be employed. In the
following sections we will compare the two possibilities from a
theoretical and numerical point of view.

\section{Relaxed numerical schemes}
\label{sec:scheme}

The numerical approximation of system \eqref{eq:sysrelphi} can be
obtained following the ideas already exploited in \cite{CNPSdegdiff}:
first obtain a semidiscrete scheme applying an IMEX time integrator
and then choose a high order non-oscillatory spatial discretization.
The reaction term will be included in the explicit part of the IMEX
scheme. This obviously is convenient only for non-stiff reaction terms.

For simplicity, here we describe the
one-dimensional case, the generalization being straightforward.

\subsection{Relaxed IMEX schemes}

We observe that system \eqref{eq:sysrelphi} is in the form
\[
 \pder{z}{t} + \pder{f(z)}{x} = g(z) + \frac{1}{\epsilon}h(z)
\]
where $z=(u,v)^T$, $f(z)=(v,\varphi^2u)^T$, $g(z)=(g(u),0)^T$
and $h(z)= (0,-v + \varepsilon\varphi^2 u_x - D A(u)u_x)^T$. When
$\varepsilon$ is small, the presence of both 
non-stiff and stiff terms, suggests the use of IMEX schemes
\cite{ARS97,KC03,PR05}. In the problems considered here, the reaction
term $g(u)$ is not stiff and hence we treat it within the explicit
portion of the IMEX scheme.

First we semi-discretize the equation in time.
Let's assume for simplicity a uniform time step $\Delta t$ and
denote with $z^n$ the numerical approximation of the variable $z$ at
time $t_n=n\Delta t$, 
for $n=0,1,\ldots$ We employ a $\nu$-stages IMEX scheme including the
$1/\varepsilon$ terms in the implicit part, giving
\begin{equation}\label{eq:IMEX}
  z^{n+1} = z^n 
     - \Delta t \sum_{i=1}^{\nu}
           \tilde{b}_i  \left[\pder{f}{x}(z^{(i)}) +g(z^{(i)})\right]
     + \frac{\Delta t}{\varepsilon} \sum_{i=1}^{\nu} b_i h(z^{(i)})
\end{equation}
where the stage values are computed as
\begin{equation}\label{eq:stage:i}
z^{(i)} = B^{(i)}
      + \frac{\Delta t}{\varepsilon} a_{i,i} h(z^{(i)})
\end{equation}
for
\begin{equation}\label{eq:stage:expl}
B^{(i)} = z^n 
        -\Delta t \sum_{k=1}^{i-1}\tilde{a}_{i,k} \left[\pder{f}{x}(z^{(k)}) +g(z^{(k)})\right]
        + \frac{\Delta t}{\varepsilon} \sum_{k=1}^{i-1} a_{i,k} h(z^{(k)})
\end{equation}
Here $({a}_{ik},{b}_i)$ and $(\tilde{a}_{ik},\tilde{b}_i)$ are a pair
of Butcher's tableaux \cite{HW1} of, respectively, a diagonally
implicit and an explicit Runge-Kutta scheme.

In this work we use the so-called relaxed schemes, that are obtained
by letting $\varepsilon\rightarrow0$ in \eqref{eq:IMEX}. The first
stage, that is \eqref{eq:stage:i} with $i=1$, is defined by
\begin{equation}\label{eq:stage:1} 
\lim_{\varepsilon\to0}
\left\{
 \left[\begin{array}{l} 
    {u^{(1)}}\\
    {v^{(1)}}\\
 \end{array}\right]
 -
 \left[\begin{array}{l} u^{n}\\v^{n}\end{array}\right]
 -
 \frac{\Delta t}{\varepsilon} a_{1,1} 
   h\left(\left[\begin{array}{l}u^{(1)}\\v^{(1)}\end{array}\right]\right)
\right\}
=0
\end{equation}
and thus
\[ {u^{(1)}=u^n} 
\qquad v^{(1)}=-A(u^{(1)})\pder{u^{(1)}}{x} .\] 
Note that, in particular, $h(z^{(1)})=0$. Now the second stage, $i=2$,
reads
\[z^{(2)} = z^n 
    - \Delta{t}\tilde{a}_{2,1}{\left[\pder{f}{x}(z^{(1)}) + g(z^{(1)})\right]} 
    + \frac{\Delta t}{\varepsilon} a_{2,1}
           \underbrace{{h(z^{(1)})}}_{{\equiv 0}}
    + \frac{\Delta t}{\varepsilon} a_{2,2} {h(z^{(2)})}.
\]
Hence the second component of $z^{(2)}$ is determined by the stiff
term of the above expression, namely $h(z^{(2)})=0$. On the other
hand, due the form of $h(z)$, there are no stiff terms in the equation
for the first component $u^{(2)}$, which is then determined by a
balance law.

Summarizing, the relaxed scheme yields an alternation of
{\em relaxation steps}
\begin{equation}\label{eq:relaxed:i}
  h(z^{(i)})=0 \qquad 
   \text{ i.e. }\, 
     v^{(i)}=-A(u^{(i)})\pder{u^{(i)}}{x}
\end{equation}
and {\em transport steps} where we advance for time $\tilde{a}_{i,k}\Delta{t}$
\begin{equation}\label{eq:transportRK}
 \pder{z}{t} + \pder{f(z)}{x} = g(z) 
\end{equation}
with initial data $z=z^{(i)}$, retain only the first component and
assign it to $u^{(i+1)}$.
Finally the value of $u^{n+1}$ is computed as \(u^n+\sum \tilde{b}_i
u^{(i)}\).

\subsection{Spatial reconstructions}

In order to have a fully discrete scheme,
we still need to specify
the space discretization. We will use discretizations based on finite
differences, in order to avoid cell coupling due to the source terms.

Recall that the IMEX technique reduces the integration to a cascade of
relaxation and transport steps. The former are the implicit parts of
\eqref{eq:stage:1} and \eqref{eq:stage:i}, while the transport steps
appear in the evaluation of the explicit terms $B^{(i)}$ in
\eqref{eq:stage:expl}. Since \eqref{eq:stage:1} and \eqref{eq:stage:i}
involve only local operations, the main task of the space discretization
is the evaluation of $\partial_x{f}$, where we will exploit the
linearity of $f$ in its arguments.

In the one-dimensional case,
 let us consider a uniform grid on $[a,b]\subset\mathbb{R}$,
\(x_j=a-\frac{h}2+jh\) for $j=1,\ldots,m$, where $h=(b-a)/m$ is the
grid spacing and $m$ the number of cells. 
We denote with $z^n_j$ the value of the quantity $z$ at time $t^n$ at
$x_j$, the centre of the $j^{\mbox{\scriptsize th}}$ computational cell. 
The fully discrete scheme
may be written as 
\[
z_j^{n+1} = z_j^n - \Delta t \sum_{i=1}^{\nu}
     \tilde{b}_i  \left(F^{(i)}_{j+1/2} - F^{(i)}_{j-1/2}\right)
    + \frac{\Delta t}{\varepsilon} \sum_{i=1}^{\nu} b_i
     g(z_j^{(i)}),  
\]
where $F^{(i)}_{j+1/2}$ are the numerical fluxes, which are the only
item that we still need to specify. It is necessary to
write the scheme in conservation form and thus, following
\cite{OS89}, we introduce the function $\hat{F}$ such that
\[
  f(z(x,t))=\frac1h \int_{x-h/2}^{x+h/2} \hat{F}(s,t) \mathrm{d}s
\]
and hence
\[
  \frac{\partial f}{\partial x}(z(x_j,t)) = \frac1h
  \left(\hat{F}(x_{j+1/2},t)-\hat{F}(x_{j-1/2},t)\right).
\]
The numerical flux function $F_{j+1/2}$ approximates
$\hat{F}(x_{j+1/2})$. 

In order to compute the numerical fluxes, for each stage value, we
reconstruct boundary extrapolated data $z^{(i)\pm}_{j+1/2}$ with a
non-oscillatory interpolation method, starting from the point values
$z^{(i)}_j$ of the 
variables at the centre of the cells. Next we apply a monotone
numerical flux to these boundary extrapolated data.

To minimize numerical viscosity we choose the Godunov flux, which, in
the present case of a linear system of equations, reduces to the upwind
flux. In order to select the upwind direction we write the linear
system with constant coefficients
\eqref{eq:transportRK} in
characteristic form. The characteristic variables relative to the
eigenvalues $\varphi,-\varphi$  are respectively
\[
      U=\frac{v+\varphi u}{2\varphi}     \qquad      
      V=\frac{\varphi{}u-v}{2\varphi}.
\]
Note that $u=U+V$. Therefore the numerical flux in characteristic
variables is
\(F_{j+1/2}=(\varphi U^-_{j+1/2},-\varphi V^+_{j+1/2})\).

The accuracy of the scheme depends on the accuracy of the
reconstruction of the boundary extrapolated data. For a first order
scheme we use a piecewise constant reconstruction such that
\(U^-_{j+1/2}=U_j\) and \(V^+_{j+1/2}=V_{j+1}\). For higher order
schemes, we use ENO or WENO reconstructions of appropriate accuracy
\cite{OS89}.
 
Since the transport steps need to be applied only to
$u^{(i)}$, we have 
\[
 u_j^{(i)}
 = u_j^n 
   -\lambda \sum_{k=1}^{i-1}\tilde{a}_{i,k}
   \left[
   \varphi \left( U^{(k)-}_{j+1/2}-U^{(k)-}_{j-1/2}
   -V^{(k)+}_{j+1/2}+V^{(k)+}_{j-1/2}\right)
   +\mathrm{\Delta t} g(u^{(k)}_j)
   \right]
\]
Finally, taking the last stage value and going back to conservative
variables, 
\[
\begin{array}{ll} 
  u_j^{n+1}=u_j^n 
  -\frac\lambda2
  \sum_{i=1}^\nu \tilde{b}_i 
  &\left([v^{(i)-}_{j+1/2}+v^{(i)+}_{j+1/2}-(v^{(i)-}_{j-1/2}+v^{(i)+}_{j-1/2})]\right.\\
  &\left.+\varphi[w^{(i)-}_{j+1/2}-w^{(i)+}_{j+1/2}-(w^{(i)-}_{j-1/2}-w^{(i)+}_{j-1/2})]
  \right)
\end{array}
\]

We wish to emphasize that the scheme reduces to the time
advancement of the single variable $u$. Although the scheme is based
on a system of three equations, the construction is used only to select
the correct upwinding for the fluxes of the relaxed scheme and the
computational cost of each time step is not affected.

\subsection{Numerical scheme}

Employing a Runge-Kutta IMEX scheme of
order $p$, can give an integration procedure for \eqref{eq:degparab}
which is of order up to $2p$ with respect to $h$ (see
\cite{CNPSdegdiff}),
because the CFL restrictions  are of parabolic type
($\mathrm{\Delta{}t}\leq{}Ch^2$). We observed that this
theoretical convergence rate can be achieved in practice, with careful
choice of the approximations of the spatial derivatives in the
discretization of \eqref{eq:relaxed:i}.

Summarizing, the relaxed schemes that we propose for the numerical
integration of \eqref{eq:degparab} consists of the following
steps. 
For each Runge-Kutta stage we need to compute the variables $v^{(i)}$
(relaxation steps). This requires the approximation of a spatial
gradient operator, for which we choose a central finite difference
operator of order at least $2p$. Close to the boundary of the domain,
there would not be enough information to compute the gradient with
high order centered finite difference formulas. We employ instead
asymmetric formulas of order $2p$.
Then we need to solve the
transport equation by diagonalizing the linear system
\eqref{eq:transportRK}, reconstructing the characteristic variables
$U^{(i)}$ and $V^{(i)}$ at cell boundaries and computing the fluxes.
Again, we must choose a spatial reconstructions of order at least $2p$
and, to avoid spurious oscillations, we employ ENO or WENO
non-oscillatory procedures. These procedures compare, for each cell,
the reconstructions obtained using different stencils and choose the
least oscillatory one (ENO) or compute a weighted linear combination
of all of them (WENO). For details, see \cite{OS89}.

Boundary conditions are enforced by extrapolating the approximate
solution and/or the stage values $u^{(\cdot)}_i$ to $1$ ghost
point located outside the domain before the computation of each
Runge-Kutta stage. This is done via a polynomial $p(x)$ of degree
$2p$, fitting the first $2p$ values of $u$ inside the domain and
satisfying the boudary condition.

\section{Numerical results}
\label{sec:tests}


\subsection{Convergence}
\begin{table}
\begin{centering}
\footnotesize
\begin{tabular}{|c|l|l|l|l|l|}
\hline
                  & N=12            & N=36            & N=108         & N=324         & N=972 \\
\hline
ENO2, RK1 &   0.00068189 &   0.00024698 & 5.9722e-005 & 1.1944e-005 & 9.4124e-007   \\
\hline
ENO3, RK2 &   0.00067142 &  0.00023703  & 1.2878e-005 & 8.8231e-007 & 3.2784e-008 \\
\hline
ENO4, RK2 &   0.00066558 &  0.00024014  & 9.7088e-006 & 2.2072e-007 & 2.0655e-009 \\
\hline
ENO5, RK3 &  0.00066357  & 0.00022958   & 4.2855e-006 & 3.3662e-008 & 1.5442e-010  \\
\hline
ENO6, RK3 &  0.0006615    & 0.0002267    & 4.2338e-006 & 1.0991e-008  & 1.833e-011   \\
\hline
WENO3, RK2 & 0.0006839  & 0.00014311  & 7.484e-005   & 1.014e-005   & 4.4262e-007 \\
\hline
WENO5, RK3 & 0.00066849 & 0.0001484   & 2.387e-005   & 2.3906e-007 & 3.7223e-010 \\
\hline
\end{tabular}

\medskip
\begin{tabular}{|c|l|l|l|l|}
\hline
	        & N=36       & N=108         & N=324         & N=972 \\
\hline
ENO2, RK1   &  0.92439     &  1.2922   &     1.465      &  2.3127  \\
\hline
ENO3, RK2   &  0.94777     &  2.6512   &    2.4401     &   2.997\\
\hline
ENO4, RK2   & 0.92793      & 2.9202    &   3.4442      &  4.2522\\
\hline
ENO5, RK3  & 0.96613       & 3.6237    &  4.4116       & 4.9011\\
\hline
ENO6, RK3  & 0.97477       & 3.6232    &   5.4194      & 5.8222\\
\hline
WENO3, RK2 & 1.4238      & 0.59004   &  1.8194       & 2.8504\\
\hline
WENO5, RK3  &  1.37       & 1.6633     &  4.1904       & 5.8847\\
\hline
\end{tabular}
\caption{Errore in norma 1 fatte sulla soluzione con 2916 punti e velocità di convergenza} 
\label{tab:errori:lin:per}
\end{centering}
\end{table}

We tested the convergence rate of our schemes with smooth initial
data. We considered equation \eqref{eq:degparab} with $A(u)=1$ and
$D=1$ and $g(u)=0$. We computed the numerical solution with the
relaxed schemes, using various combinations of spatial reconstructions
and time integrators, starting from the initial data $u(x)=\sin(2\pi
x)$ for $x\in[0,1]$. The errors and convergence rates are reported in
Table \ref{tab:errori:lin:per}. This test shows that the schemes can
reach the predicted order of convergence.

\subsection{Travelling waves}

\begin{figure}
\includegraphics[width=0.45\textwidth]{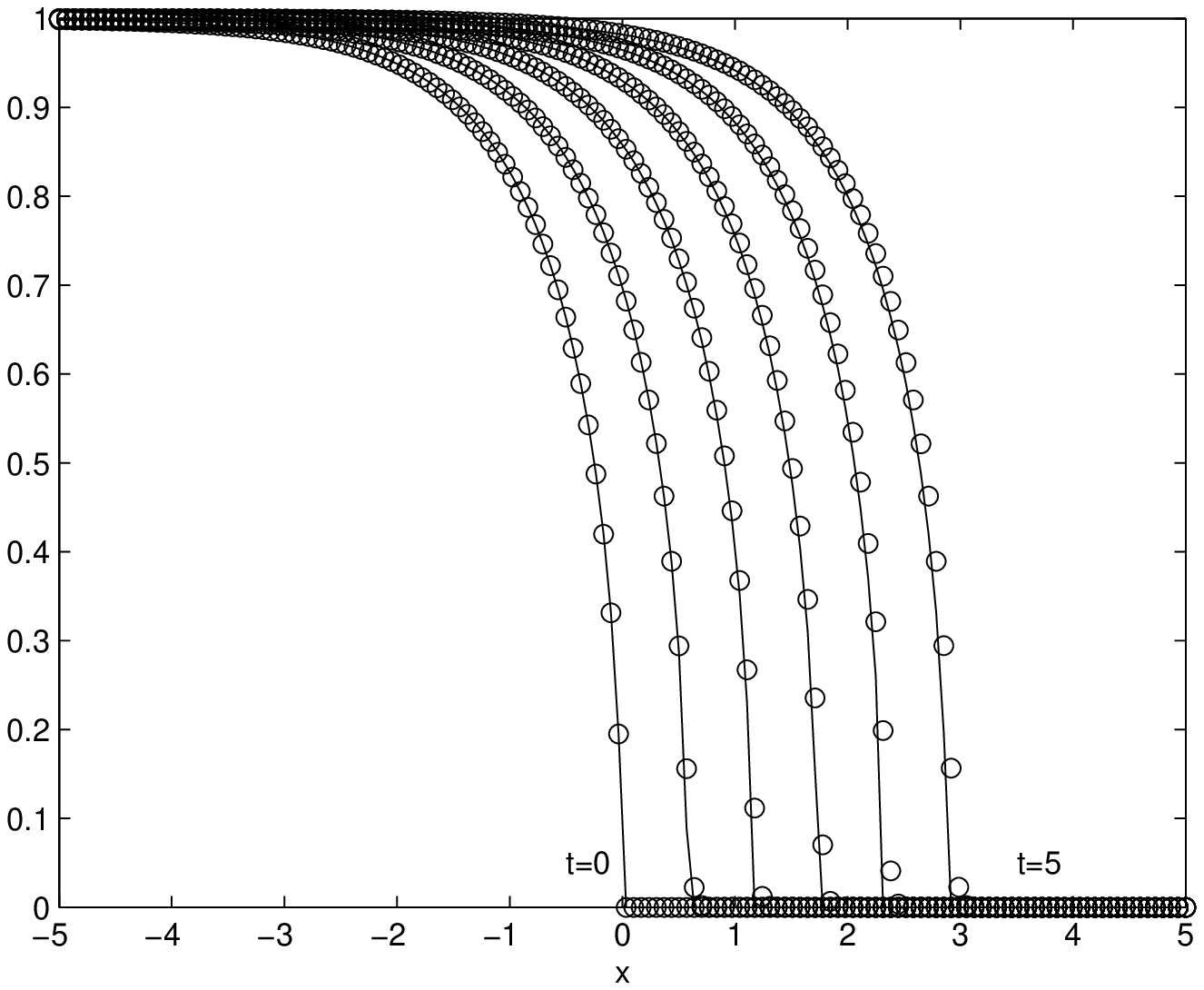}
\hfill
\includegraphics[width=0.45\textwidth]{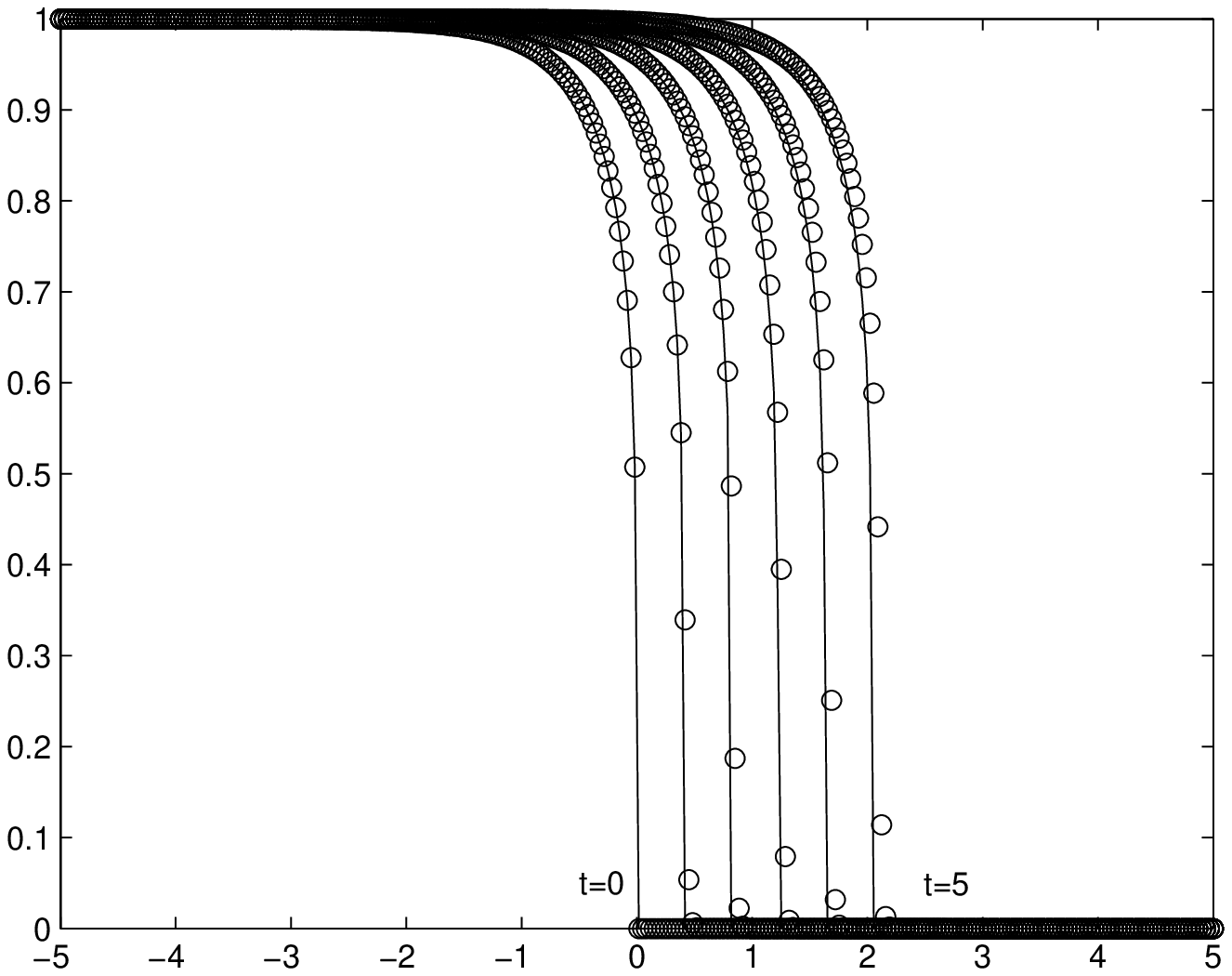}

\caption{Travelling waves with $\alpha=2$ (left) and $\alpha=5$
  (right). The circles represent the numerical solution at
  $t=0,1,2,3,4,5$ and the solid lines are the corresponding exact
  solutions \eqref{eq:genfk:tw}.}
\label{fig:tw}
\end{figure}

We consider the following generalization of the Fisher-Kolmogoroff
equation \eqref{eq:genfk}
\begin{equation} \label{eq:genfk}
 \pder{u}{t} = u^p(1-u^q) +  \pder{}{x}\left(u^m\pder{u}{x}\right) 
\end{equation}

The existence of travelling waves can be proved for a wide range of
parameters $p,q,m$ \cite{Murray1}. For example in the paper
\cite{Newman83}
gives the expression of such a wave for the case $p=1$
and $q=m=\alpha$. The speed $c$ and the exact solution are:
\begin{equation}\label{eq:genfk:tw}
  c = \frac1{\sqrt{1+\alpha}}
  \qquad
  u(t,x) = 
  \left[
    \left(
      1 - e^{\frac{\alpha(x-ct)}{\sqrt{1+\alpha}}}
    \right)^{\frac1{\alpha}}
  \right]_+
\end{equation}

Figure \ref{fig:tw} shows the numerical solutions obtained for $t=5$,
using $300$ points in the interval $x\in[-5,5]$. The initial data and
the exact solutions shown are defined according to
\eqref{eq:genfk:tw}. The accuracy of the schemes is good in both
cases. In the case with $\alpha=5$ there is a lower accuracy around
the corner, due to the higher degeneracy of the diffusion term and the
more pronounced stiffness of the reaction term.

Table \ref{tab:tw} shows the errors and convergence rates, measured at
$t=5$, comparing the results of the schemes presented in this paper
and those obtained with the schemes described in
\cite{CNPSdegdiff}. In order to apply the schemes of
\cite{CNPSdegdiff}, equation \eqref{eq:genfk} was rewritten in
the ``conservative form''
$u_t=[u^m/m]_{xx}+u(1-u^\alpha)$, where $m=\alpha+1$.
We point out that, although we employed the Heun scheme and ENO
reconstruztions of order $3$, the convergence rate is limited by the
$C^{(0)}$ regularity of the solution. 
The comparison shows that the schemes for the conservative form
perform slightly better than those described here for the
nonconservative form, both regarding the convergence rate and the
absolute value of the errors. Both schemes could be improved by
integrating the reaction term with the diagonnally implicit part of
the IMEX scheme, expecially in the test with $\alpha=5$.

Thus, whenever the equation can be put in the conservative form
required by the schemes of \cite{CNPSdegdiff}, those are to be
preferred. In the other cases, which frequently occour in more than
one dimensions or in the case of systems of equations, the present
scheme is a good choice for the numerical integration of the
differential problem.

\begin{table}
\begin{tabular}{|c|c|c|c|c|c|c|c|c|}
\hline
 & \multicolumn{4}{c|}{$\alpha=2$} & \multicolumn{4}{c|}{$\alpha=5$} \\
\hline
    & \multicolumn{2}{c|}{Section \ref{sec:scheme}} &  \multicolumn{2}{c|}{\cite{CNPSdegdiff}}
    & \multicolumn{2}{c|}{Section \ref{sec:scheme}} &  \multicolumn{2}{c|}{\cite{CNPSdegdiff}} \\
\hline
$N$   & $\|E\|_2$ & rate & $\|E\|_2$ & rate& $\|E\|_2$ & rate& $\|E\|_2$ & rate\\
\hline
30 & 2.3515e-01 & & 6.9686e-02 & & 4.5251e-01 & &  5.0452e-02 & \\
\hline
60 & 6.9899e-02 & -1.75& 4.7352e-03 & -3.87 & 2.0266e-01 & -1.15 & 2.7865e-02 & -0.85\\
\hline
120& 2.4679e-02 & -1.50& 4.1303e-03 & -0.19 & 9.0192e-02 & -1.16 & 1.3893e-02 & -1.00\\
\hline
240& 6.0385e-03 & -2.03& 7.0980e-04 & -2.54 & 3.2529e-02 & -1.47 & 3.1202e-03 & -2.15\\
\hline
480& 2.1158e-03 & -1.51& 3.2753e-04 & -1.11 & 1.2977e-02 & -1.32 & 1.1481e-03 & -1.44\\
\hline
avg. rate && -1.71 && -1.82 && -1.28 && -1.40\\
\hline
\end{tabular}
\caption{Errors and convergence rates on the travelling wave solution
  \eqref{eq:genfk:tw}. For $\alpha=2$ and $\alpha=5$ we compare the
  results of the schemes proposed in this paper with those described
  in \cite{CNPSdegdiff}. The last row is the average convergence rate.}
\label{tab:tw}
\end{table}

\subsection{Two dimensional tests}
\begin{figure}
\includegraphics[angle=-90,width=.95\textwidth]{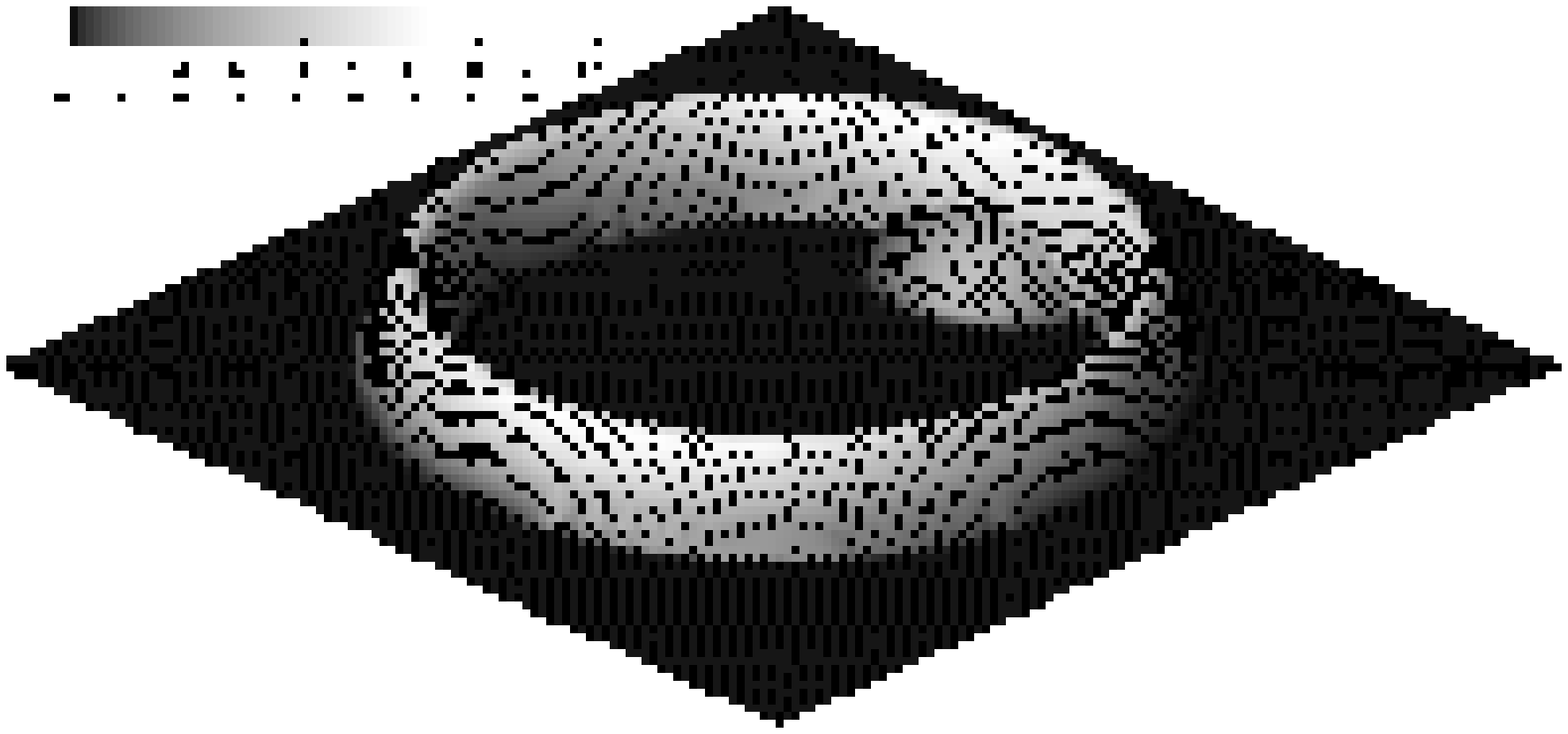}

\includegraphics[width=.24\textwidth]{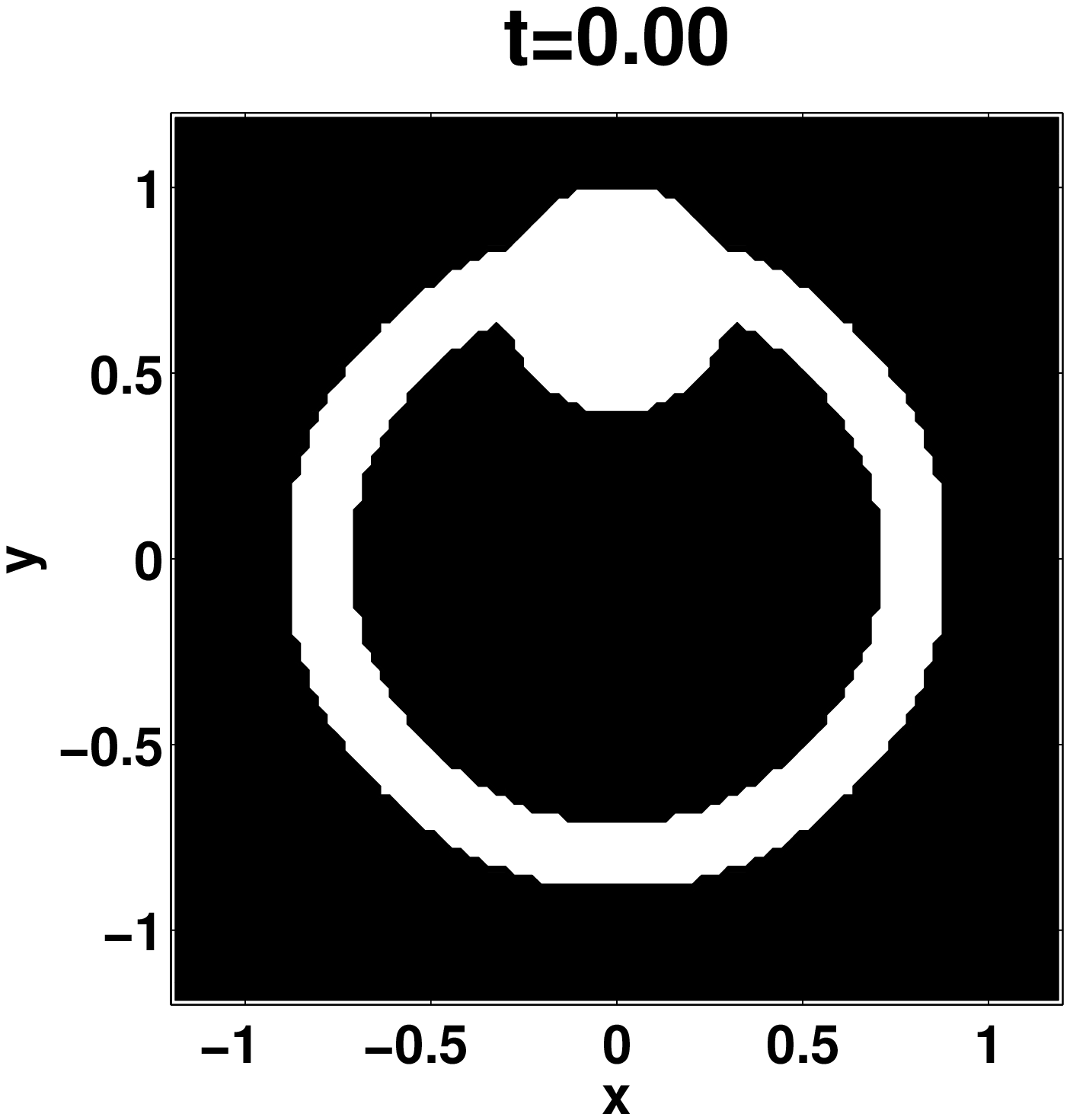}
\includegraphics[width=.24\textwidth]{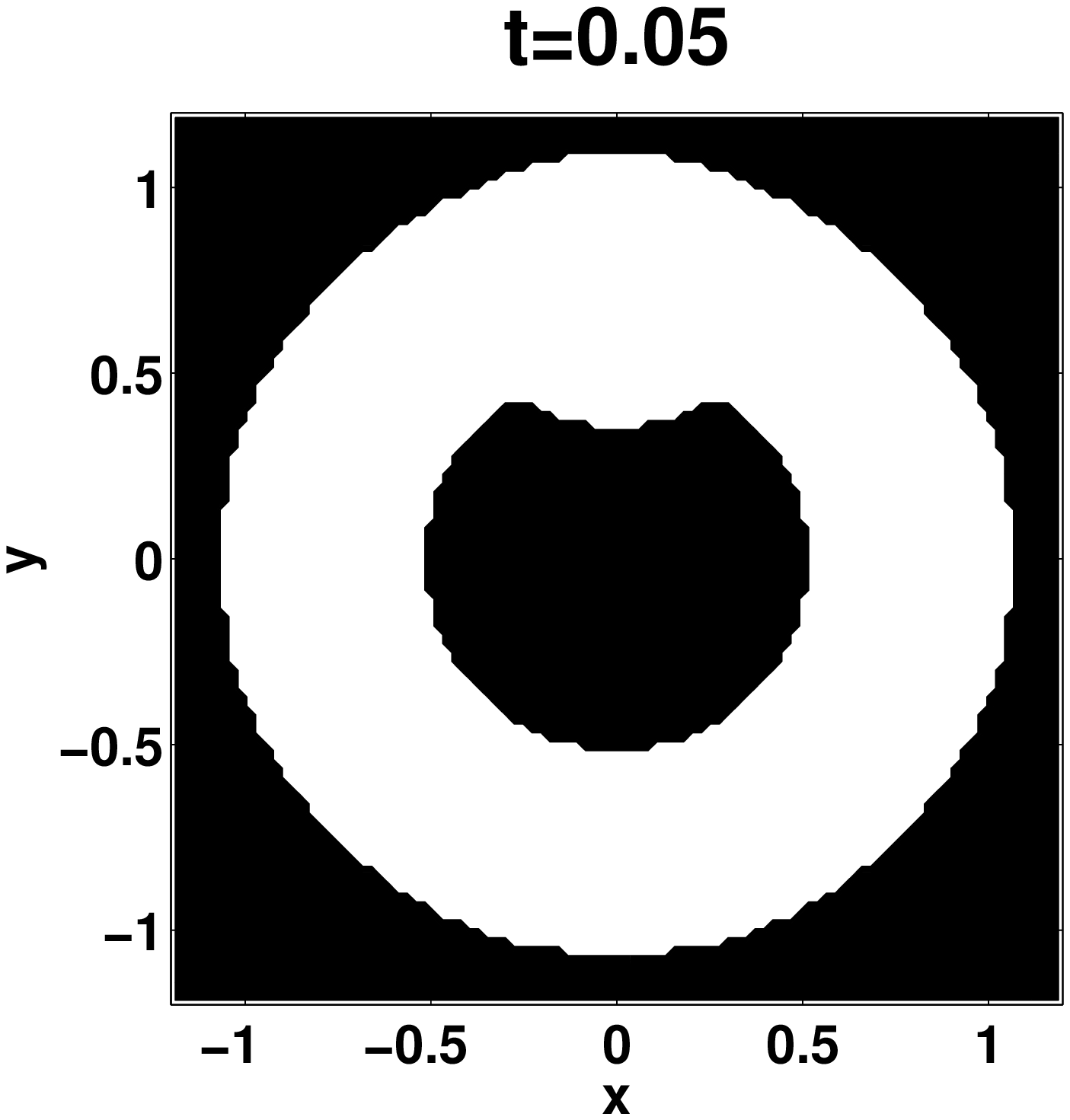}
\includegraphics[width=.24\textwidth]{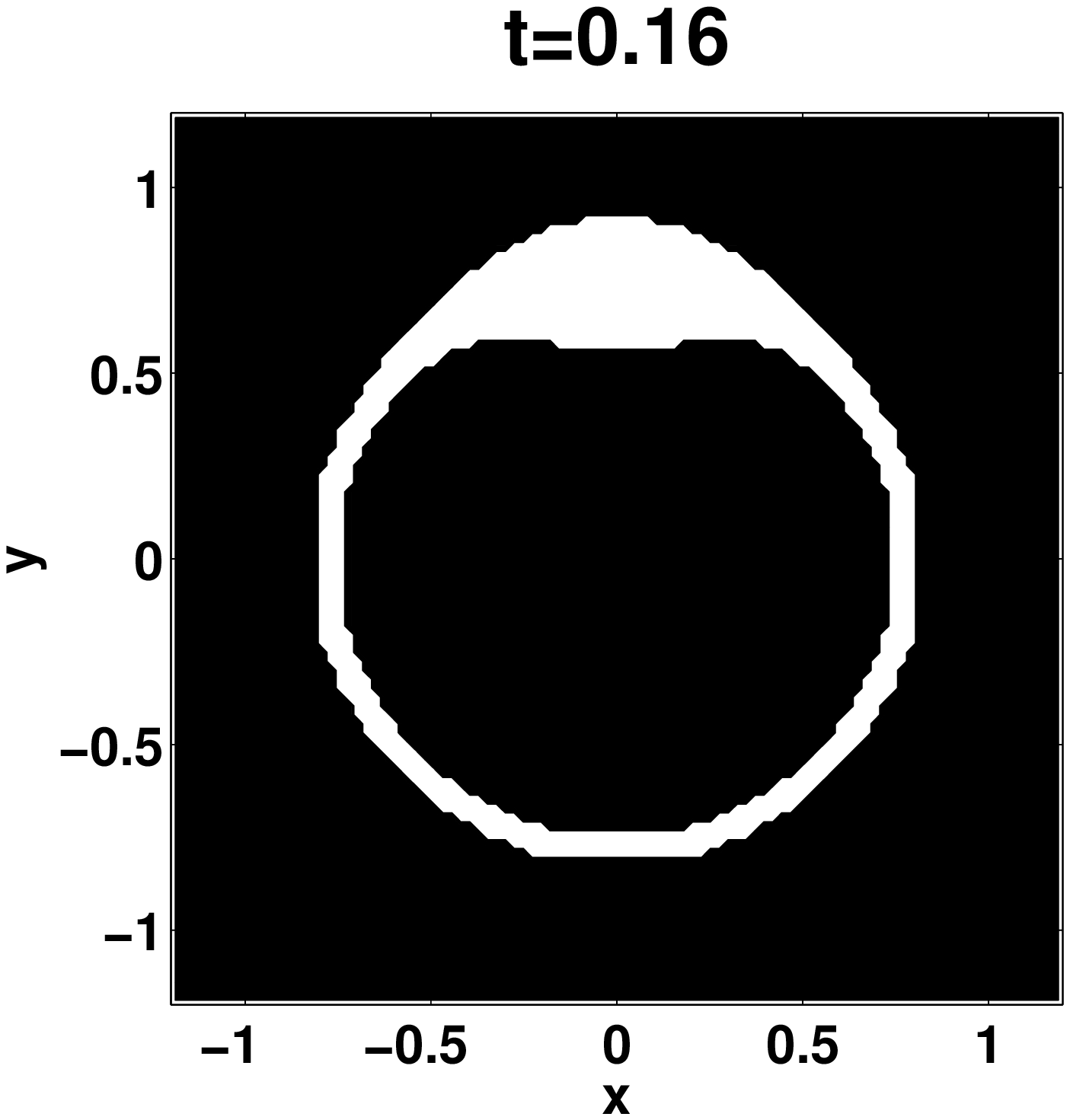}
\includegraphics[width=.24\textwidth]{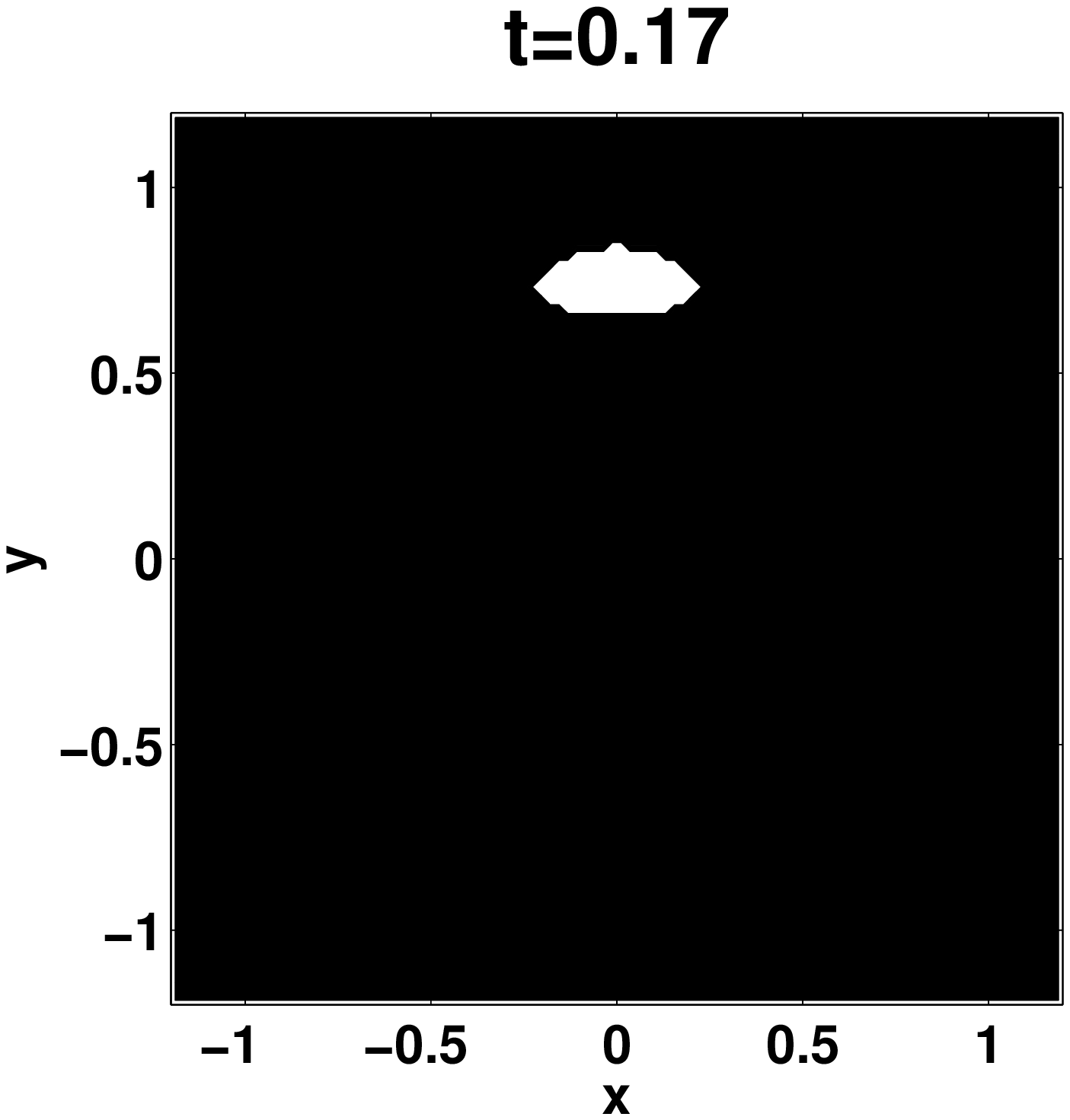}

\caption{Initial data (top) and supports (bottom) of the numerical
  solution of \eqref{eq:mik} at different times.}
\label{fig:testanellobolla}
\end{figure}

We tested the multi-dimensional version of our integration scheme on
the equation
\begin{equation} \label{eq:mik}
 \pder{u}{t} = \mathrm{\Delta}(u^m) -cu^p 
\end{equation}
for $x\in[-2,2]^2$ and $t\geq0$. The above equation presents
interesting finite-time extinction phenomena, reported in
\cite{Mik95}.

We tested the ability to reproduce the extinction phenomenon and the
persistence of asymmetry in the initial datum along the evolution. We
took $u(x,y,0)$ as a radially symmetric function with a small
perturbation, see Figure \ref{fig:testanellobolla}, and evolved it
with the same equation and parameters as before, until extinction.
Figure \ref{fig:testanellobolla} shows clearly that the initial
perturbation of the radial symmetry is maintained until the solution
vanishes.

\subsection{Application to a system of PDEs}
\begin{figure}
\begin{tabular}{cc}
\includegraphics[width=.45\textwidth]{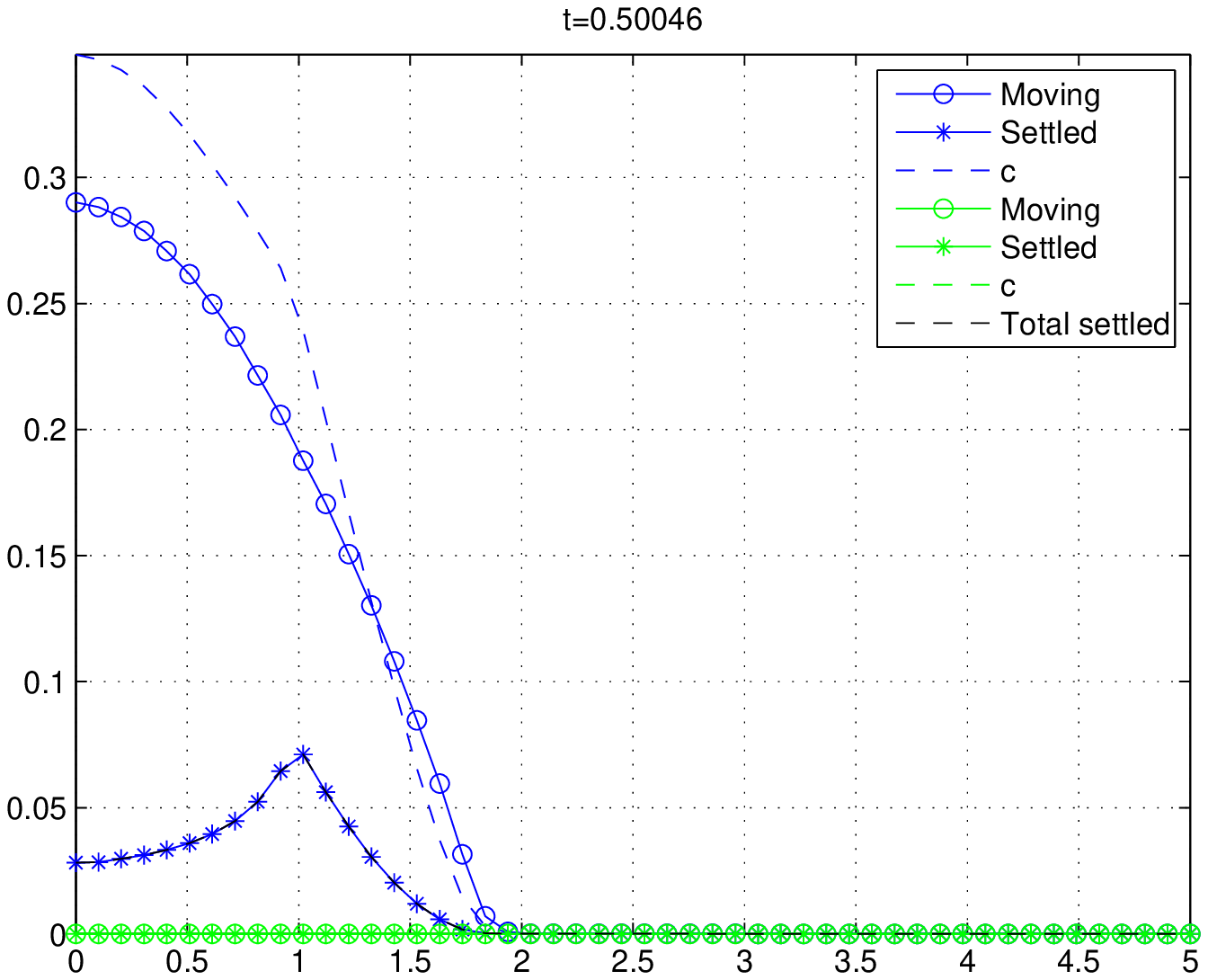}&
\includegraphics[width=.45\textwidth]{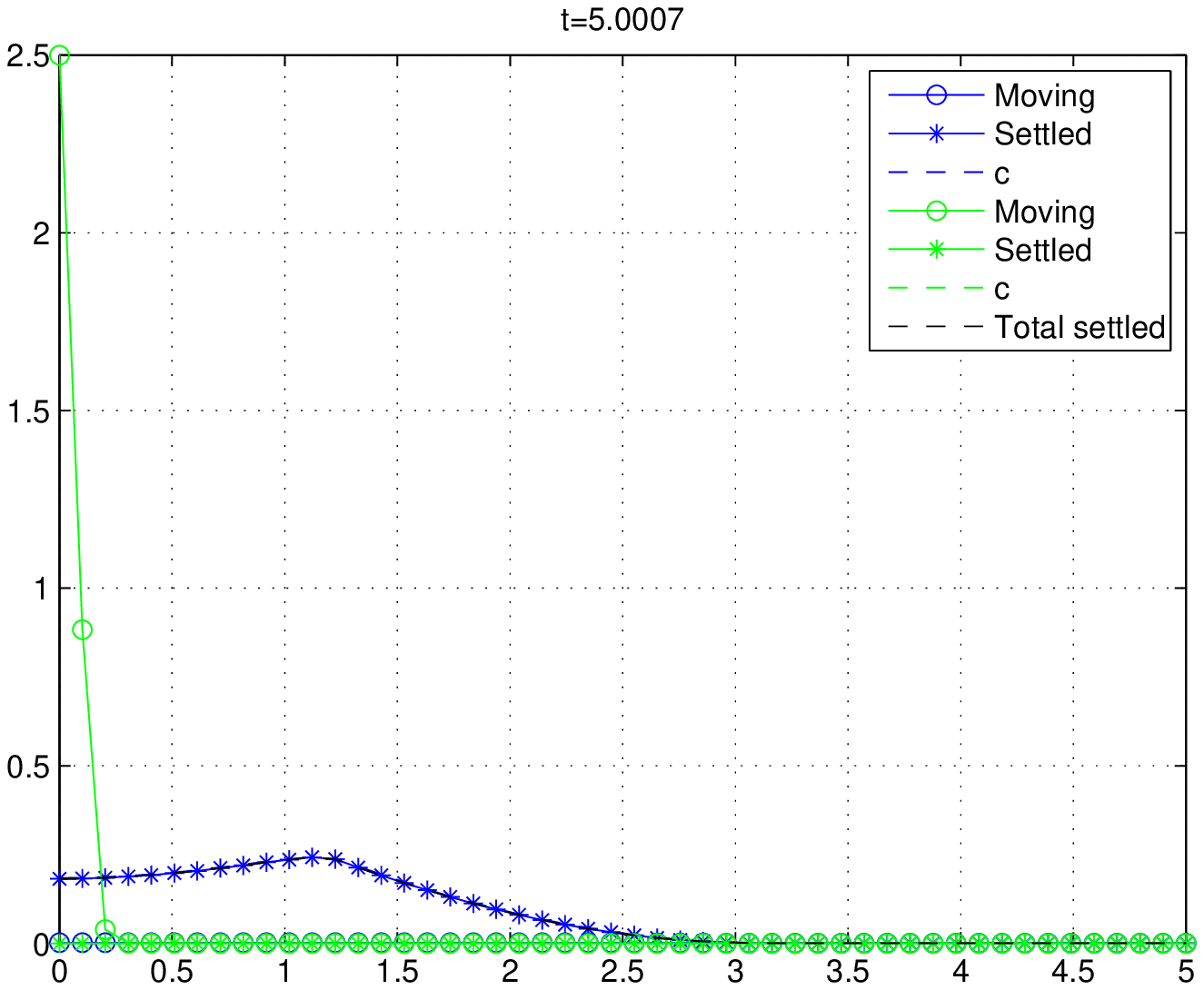}\\
\includegraphics[width=.45\textwidth]{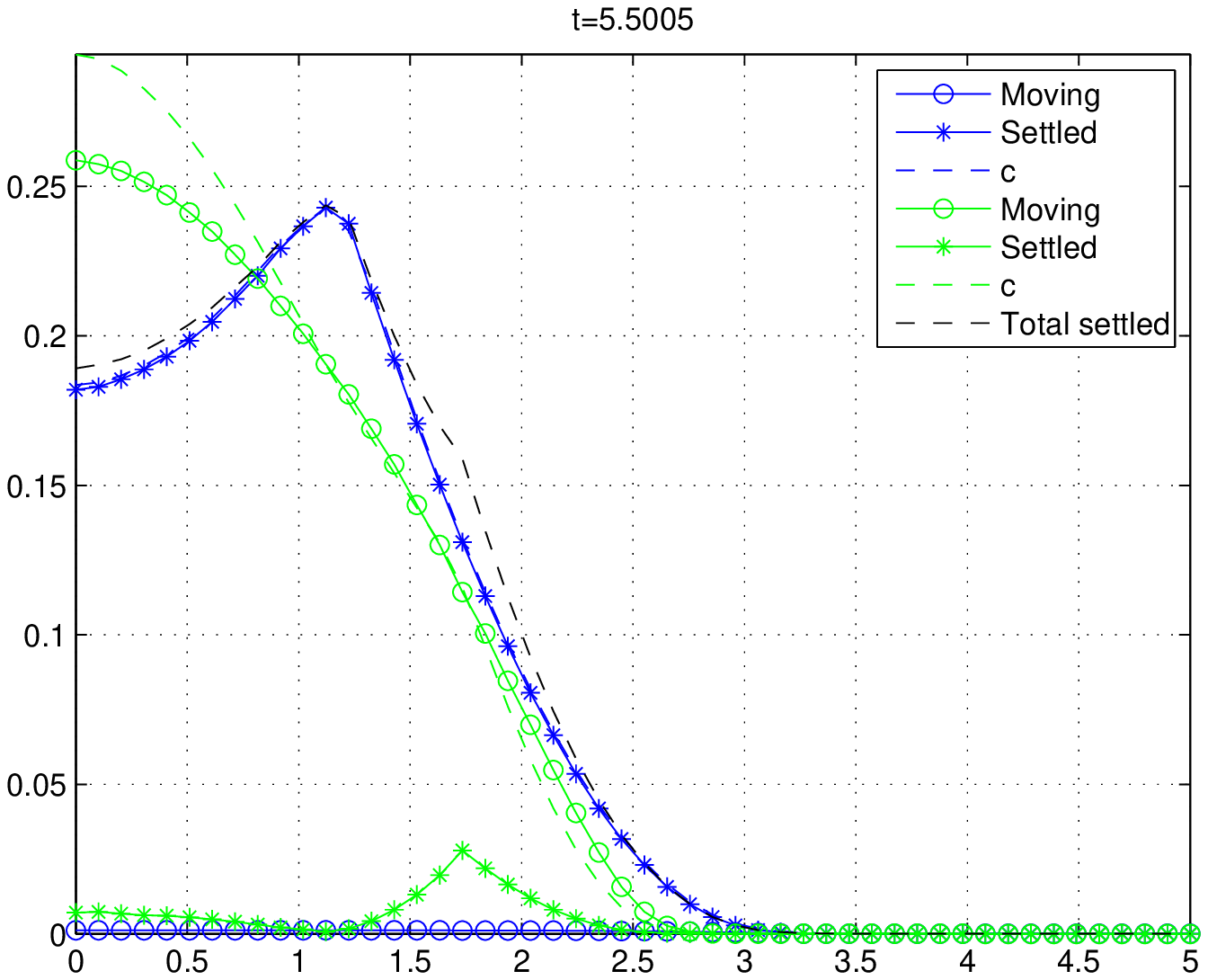}&
\includegraphics[width=.45\textwidth]{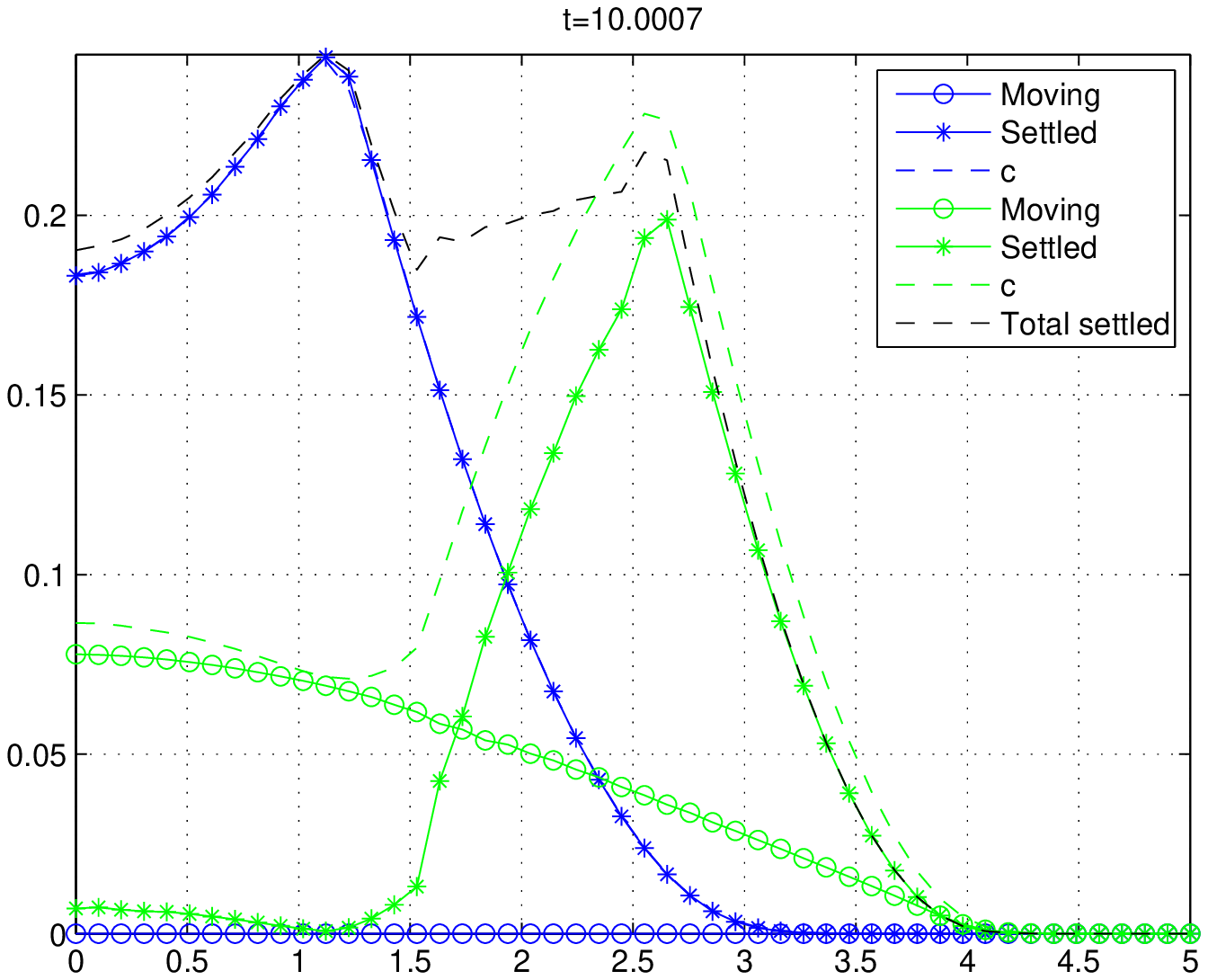}\\
\includegraphics[width=.45\textwidth]{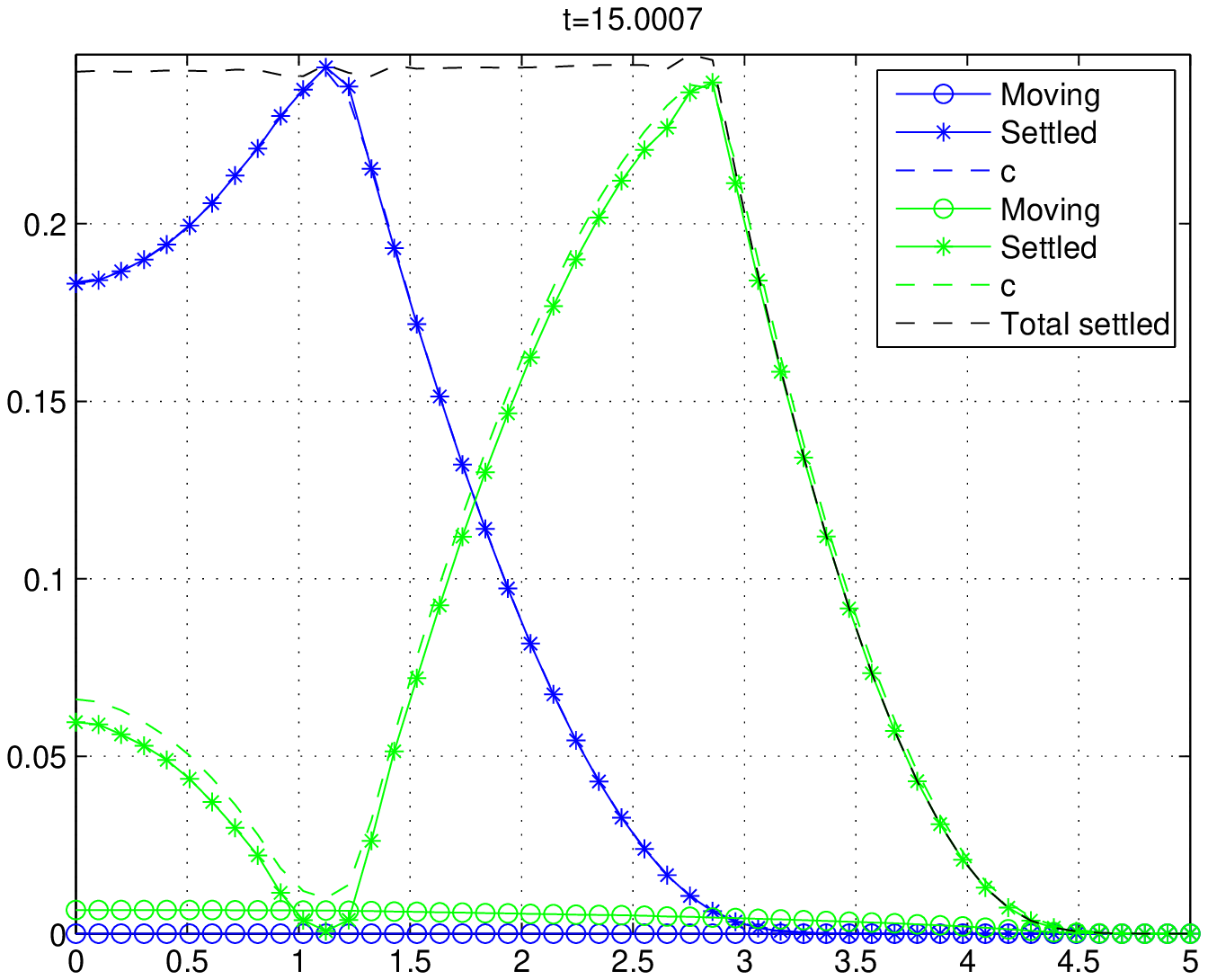}&
\includegraphics[width=.45\textwidth]{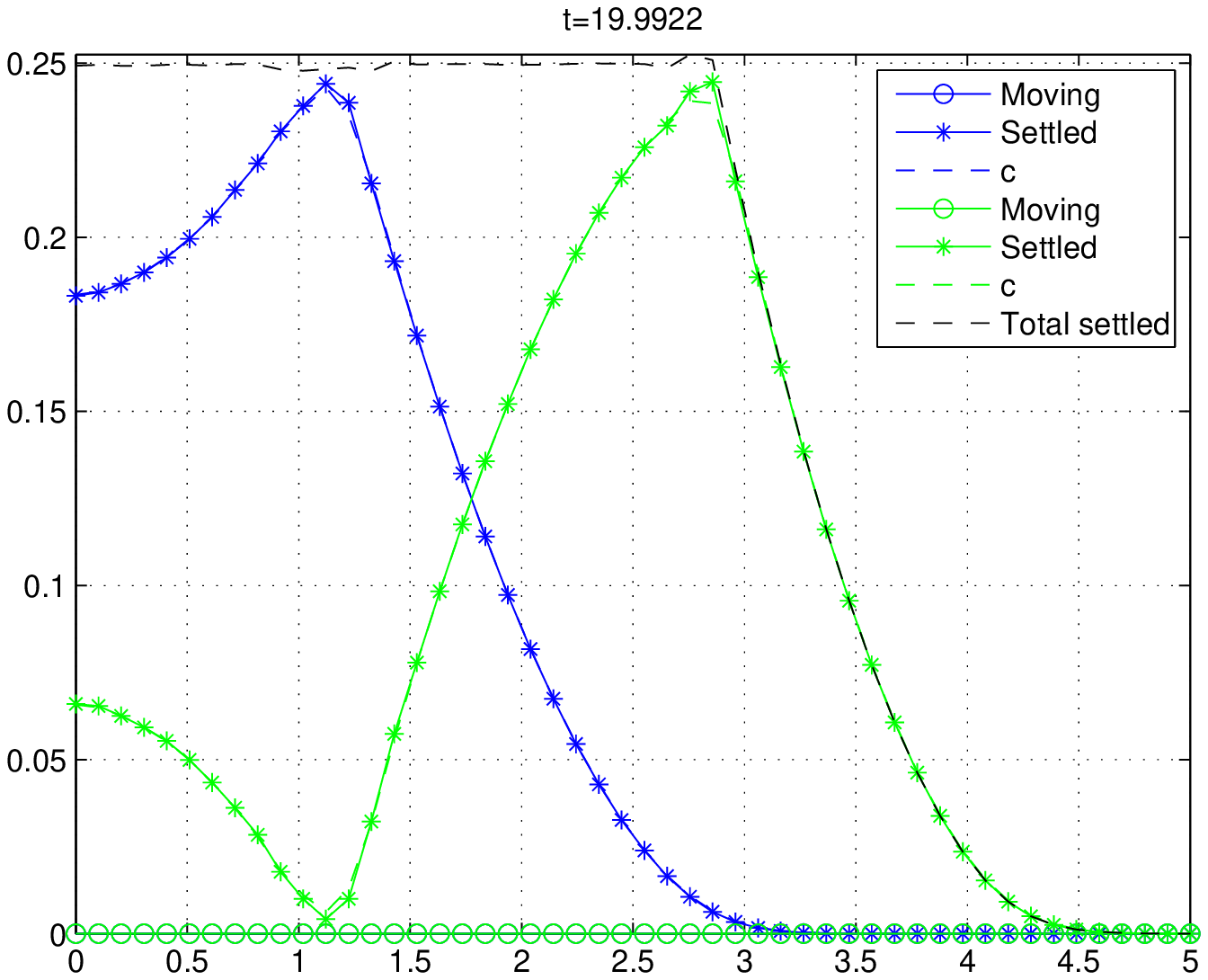}
\end{tabular}
\caption{Model \eqref{Frogs:model} for frog populations in
  \cite{TLB07}. Solutions at $t=0.5$, $t=5$, $t=5.5$, $t=10$, $t=15$
  and $t=20$. Note the different scale in the top right panel. In
  blue: the $u$ population released at $x=0$ at $t=0$; in green: the
  $v$ population, that is released at $x=0$ after the first has
  settled ($t=5$ here). The black dashed line is the total settled
  population, i.e. $u_s+v_s$. 
}
\label{fig:frogs}
\end{figure}

We consider now an ecological model described in \cite{TLB07}. It
consists in a system of PDEs describing the dispersal and settling of
frogs. The model splits the population of frogs in two classes: the
ones that have already settled and those that are still migrating.
The moving class diffuses with rate $D$ and settles at rate $S$. Once
an individual is in the settled class, they cannot migrate any more.
The authors of \cite{TLB07} study various mechanisms of dispersal and
settling: we describe here only those that they think relevant in the
case study considered in the paper, that is density-dependent
dispersal and constant or $u$-dependent settling rate.

\begin{subequations} \label{Frogs:model} 
  Let $u_m(t,x)$ and $u_s(t,x)$ be the population density of migrating
  and settled animals, $D(u)$ and $S(u)$ suitable functions describing
  the dispersal and settling mechanism. The population is described by
  the following parabolic equations.
\begin{align}
\label{Frogs:um}
\pder{u_m}{t} &= \mu \div\left(D(u_m+u_s)\grad(u_m)\right) 
                    -S(u_m+u_s)u_m\\
\label{Frogs:us}
\pder{u_s}{t} &= S(u_m+u_s)u_m
\end{align}
Initial conditions mimic a population concentrated around $x=0$ for
$u_m(0,x)$ and $u_s(0,x)=0$.  In the case $D(u)=u/u_0$ (for constant
$u_0$) considered in the paper, the diffusion equation
\eqref{Frogs:um} is degenerate and gives rise to Barenblatt-like
profiles for the density $u_m$ of the moving population. We point out
that a very accurate numerical approximation is very important, since
there is no mechanism in the model to move the population once it has
settled.

If a second population $v$ is released in the same spot, the model is
augmented with the following equations: 
\begin{align}
\label{Frogs:cu}
\pder{c_u}{t} &= \alpha\Laplacian_x(c_u)+\beta(u_m+u_s-c_u)\\
\label{Frogs:vm}
\pder{v_m}{t} &= \mu \div\left(D(v_m+v_s)\grad(v_m)\right) 
                 +\gamma \div\left(v_m\grad(c_u)\right)
                 -S(u_s+v_m+v_s)v_m\\
\label{Frogs:vs}
\pder{v_s}{t} &= S(u_s+v_m+v_s)u_m
\end{align}
Here  $v_m$ and $v_s$ are the densities of  migrating and settled
animals belonging to the second population, while $c_u$ is a
diffusible pheromone released by the first population and discouraging
the $v$ population from staying in the area where the $u$ density is
higher. 
\end{subequations}

The case study of \cite{TLB07} is about a release of a batch of frogs
($v$) at the same spot where an earlier release ($u$) had taken place,
at a time when the frogs from the first batch have already settled
down.  In order to reproduce the situation of the case study of
\cite{TLB07}, one has to evolve in time the system \eqref{Frogs:model}
starting with an initial nonzero $u_m$ population while
$u_s=v_m=v_s=0$ and adding a nonzero $v_m$ density after the $u$
population had time to settle settled (that is when $u_m$ is
negligible).

We employed free-flow boundary conditions and chose parameters as
suggested in \cite{TLB06}:
\[
\mu = 1 \qquad 
\gamma = 0 \qquad
\alpha = 0.01 \qquad
\beta  = 10 \qquad
\]
The initial population is taken to be
\[u_m(0,x)= 2.5 e^{-100x^2}\]
and we chose
\[ D(u) = \frac{u}{0.25}
   \qquad
   S(u,c)= \chi\left(1-\frac{u}{0.25}\right)
\]
where $\chi(x)=1$ when $x>0$ and $\chi(x)=0$ otherwise. 

The $u$ population settles in the time lapse from $t=0$ to $t=5$ with
maximum density located some distance apart from the original release
spot ($x=0$ in the simulation). After this, the release of the second
batch is simulated setting, at $t=5$, $v_m(5,x)=u_m(0,x)$.  Thus the
same initial condition is employed for both batches.
    
Figure \ref{fig:frogs} clearly shows that this second population first
moves further away than the peak of the $u_s$ distribution and start
settling there. Only afterwards the $v$ population settles also near
$x=0$.

The panel in figure \ref{fig:frogs} shows the evolution of the system,
from the early $u$ dispersal that gives rise to a Barenblatt-like
profile for $u_m$ (upper left), to the release of the $v$ population
when $u_m\simeq0$ (upper right), the early migration and settling of
the second population around $x=2$ (middle), the beginning of the
secondary settling at $x=0$ (bottom left) and the final steady state
(bottom right).

We reproduced with our method also all the other cases studied in
\cite{TLB06}. We point out that the standard integration procedure
used for the simulations in \cite{TLB06} required a smoothing of the
$\chi$ function and, despite that, produces unrealistic wiggles in the
population densities, that are not seen with the method presented here.

\section{Conclusions}

We have proposed and analyzed relaxed schemes for nonlinear degenerate
reaction diffusion equations.  By using suitable discretization in
space and time, namely ENO/WENO non-oscillatory reconstructions for
numerical fluxes and IMEX Runge-Kutta schemes for time integration, we
have obtained a class of high order schemes. The theoretical
convergence analysis for the semidiscrete scheme and the stability for
the fully discrete schemes have been studied by us, for the case of
nonlinear diffusion, in \cite{CNPSdegdiff}.

Here we tested these schemes on travelling waves solutions and in
cases where the solution vanishes in finite time. In all cases we
observed a very good agreement with known properties of the exact
solutions.

Moreover we have shown how the accurate approximation of the fronts
arising in degenerate parabolic equations can lead to improved
accuracy in the simulation of an ecologically relevant model.

\bibliographystyle{alpha}
\bibliography{zennoncons}

\end{document}